\def\x{{\bf{x}}}
\def\X{{\bf{X}}}

\def\pr{{\hbox{Pr}}}

\def\qed{\hfil\hfil{\vrule height7pt width6pt depth0pt}\hfil}

\documentclass[12pt]{article}
\newtheorem{theorem}{Theorem}[section]

\newtheorem{proposition}[theorem]{Proposition}

\newtheorem{algorithm}[theorem]{Algorithm}

\usepackage{amssymb}
\usepackage{amsmath}
\usepackage{graphicx}% Include figure files in eps format.
\usepackage{bm}% bold math.
\usepackage{epstopdf}

\title{Simulation with Fluctuating and Singular Rates}

\author{ 
Farzin Barekat\thanks{Mathematics Department, University of California at Los Angeles,
Los Angeles, CA 90095-1555 USA. fbarekat@math.ucla.edu. Research supported by DOE grant DE-FG02-05ER25710.} ~
and  
Russel Caflisch\thanks{Mathematics Department, University of California at Los Angeles,
Los Angeles, CA 90095-1555 USA. caflisch@math.ucla.edu. Research supported by DOE grant DE-FG02-05ER25710.}
}

\date{\ }

\begin{document}
\maketitle

\begin{abstract}
In this paper we present a method to generate independent samples for a general random variable, either continuous or discrete. The algorithm is an extension of the acceptance-rejection method, and it is particularly useful for kinetic simulation in which the rates are fluctuating in time and have singular limits, as occurs for example in simulation of recombination interactions in a plasma. Although it depends on some additional requirements, the new method  is  easy to implement and rejects less samples than the acceptance-rejection method. 
\end{abstract}

\section{Introduction}\label{section:Introduction}

\begin{figure} [!htp]
    \centering
        \includegraphics[width=16cm]{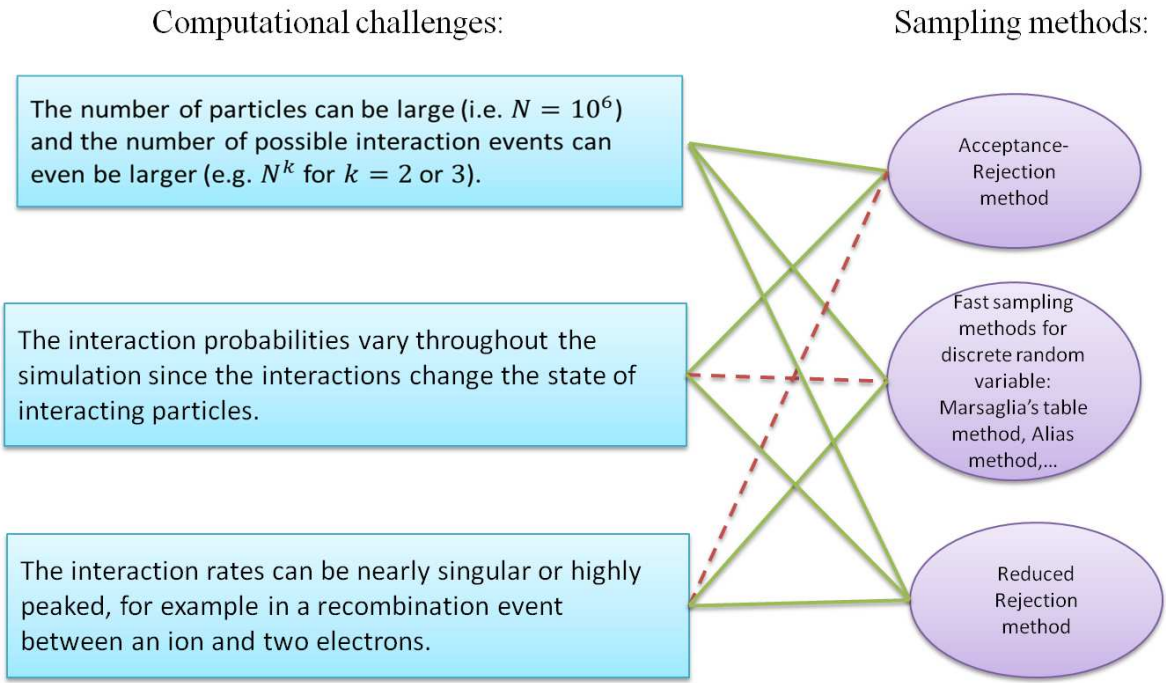} 
    \caption{This figure illustrates the computational challenges involved in sampling interactions of numerical particles, and how different methods can handle them. Broken line represents challenges for which the method becomes computationally inefficient, whereas, the solid line represents the challenges for which the method is still computationally efficient.}
    \label{figure:challenge}
\end{figure}

Kinetic transport for a gas or plasma involves particle interactions such as collisions, excitation/deexcitation and ionization/recombination. Simulation of these interactions is most often performed using the Direct Simulation Monte Carlo (DSMC) method \cite{BIB:Bird1994} or one of its variants, in which the actual particle distribution   is represented by a relatively small number of numerical particles, each of which is characterized by state variables, such as position $x$ and energy $E$. Interactions between the numerical particles are performed by random selection of  the interacting particles and the interaction parameters, depending on the interaction rates. Correctly sampling these interactions involves several computational challenges: First the number $N$ of particles can be large (e.g., $N=10^6$) and the number of possible interaction events can be even larger (e.g., $N^k$ for $k=2$ or $3$). Second, the interaction probabilities vary throughout the simulation since  interactions change the state of the interacting particles. These two difficulties are routinely overcome using acceptance-rejection sampling. Third, the interaction rates can be nearly singular, for example in a recombination event between an ion and two electrons (described in more detail in Section \ref{section:Example2}). This creates a wide range of interaction rates that makes acceptance-rejection computationally intractable. Figure \ref{figure:challenge} illustrates these challenges and how different methods can handle them.
The sampling method presented here, which we call Reduced Rejection, was developed to overcome the challenges of a large number of interaction events with fluctuating and singular rates.

%The development of the sampling method in this paper, called PQ-sampling, was motivated by a problem in plasma physics described in section \ref{section:Example2}. The existence of singularity and time-varying characteristic of the probability distribution function in this physics problem, as well as, the requirement to generate uncorrelated samples renders the use of other conventional sampling methods inefficient. The PQ-sampling method is devised to overcome the challenges in generating uncorrelated samples for random variables that have probability distribution function with aforementioned properties.

%Simulating a random variable is required for many practical purposes and has applications to variety of fields. As a result, it is important to have algorithms that simulate a random variable efficiently. 

Simulation of kinetics requires sampling methods that generate independent samples. This rules out  Markov Chain Monte Carlo  schemes, such as Metropolis--Hastings, Gibbs sampling, and Slice sampling.
 Although these methods  are very powerful and are used very often, this paper focuses on  sampling methods that generate independent samples.

%The sampling methods that generate independent samples can be further divided into those methods used for discrete random variables and methods applicable to continuous random variables\footnote{General speaking, because the density function of a discrete random variable can be thought as a special continuous density function, any method used for continuous random variable can be used for discrete random variable, however, it might not be as efficient for discrete random variables as other methods.}.
There are several efficient algorithms for simulation of discrete random variables, notably Marsaglia's table method \cite{BIB:TableMethod} and the Alias method \cite{BIB:AliasMethod,BIB:Walker}. 
%and Uniform Ratio method \cite{BIB:UniformRatio}. 
However, these methods require pre-processing time and, therefore, are not efficient for sampling from a random variable whose probability  function changes during the simulation. For continuous random variables there are several different algorithms; nevertheless, each of these algorithms has its own constraints. For example, Inverse Transform Sampling method requires knowledge of the cumulative distribution function and evaluation of its inverse, Box-Muller  only applies to a normal distribution, and Ziggurat algorithm \cite{BIB:Ziggurat} can be used for random variables that have monotone decreasing (or symmetric unimodal) density function.

The algorithm of choice for general (both continuous and discrete) random variables that generates independent samples and does not require preprocessing time is acceptance-rejection method (see for example \cite{BIB:Caflisch}). Let $q(x)$ be a real-valued function on the sample space. Let $I[q]$ denote the expectation of function $q(x)$. By sampling according to function $q(x)$ we mean to sample using the probability distribution function $q(x)/I[q]$. We say function $q(x)$ \emph{encloses} function $p(x)$ if $p(x)\leq q(x)$ for all $x$ in the sample space. The idea of acceptance-rejection method is to find a proposal function $q(x)$ that encloses function $p(x)$. Suppose we already have a mechanism to sample according to $q(x)$, then acceptance-rejection algorithm enables us to sample according to $p(x)$. In most cases the constant function is used as the proposal function $q(x)$. The main drawback of acceptance-rejection method is that it might  reject many samples. Indeed the ratio of the number of rejected samples to the number of accepted samples is approximately equal to the ratio of the area between curves $q(x)$ and $p(x)$ to the area under the curve $p(x)$.

For many given distributions, finding a good proposal function that encloses it without leading to many rejected samples is difficult. One extension to acceptance-rejection method is Adaptive Rejection Sampling \cite{BIB:GilksWild}. The basic idea of Adaptive Rejection Sampling is to construct proposal function $q(x)$ that encloses the given distribution by concatenating segments of one or more exponential distributions. As the algorithm proceeds, it successively updates the proposal function $q(x)$ to correspond more closely to the given distribution. %Other generalization of acceptance-rejection are described in section \ref{section:Comparison}.
Another extension to Acceptance-rejection method is economical method \cite{BIB:Deak}. This method is basically a generalization of Alias method for continuous distributions. In this method, one needs to define a specific transformation that maps $\{x:p(x)>q(x)\}$ to $\{x:p(x)\leq q(x)\}$. Although this method produces no rejection, finding the required transformation is difficult in general.

%One requirement of adaptive rejection sampling is that the function $p(x)$ that we want to sample according to must be log-concave. Adaptive rejection Metropolis sampling \cite{BIB:GilksBestTan} is a generalization of adaptive rejection sampling that has a Metropolis step to accommodate distributions that are not log-concave.\\

In the Reduced Rejection  method we sample according to a given function $p(x)$ based on a proposal function $q(x)$. In contrast to the acceptance-rejection method, Reduced Rejection sampling does not require $q(x)$ to enclose $p(x)$ (i.e. it allows $p(x)> q(x)$ for some $x$). On the other hand, Reduced Rejection sampling requires some extra knowledge about the functions $p(x)$ and $q(x)$.% Namely, we need to know the values of $I[p]$, $I[q]$ and have a mechanism to sample according to $p(x)-q(x)$ from those $x$ for which $p(x)>q(x)$. \\

The Reduced Rejection sampling method can be applied to a wide range of sampling problems (for both continuous and discrete random variables) and in many examples is more efficient than customary methods (three examples are provided in sections \ref{section:Example1},  \ref{section:Example2} and \ref{section:Example3}). In particular, Reduced Rejection sampling requires no pre-processing time and consequently is  suitable for simulations in which $p(x)$ is changing constantly  (see section \ref{section:Comparison} for an elaboration on this point and sections \ref{section:Example2} and \ref{section:Example3} for examples of simulations with fluctuating $p(x)$). Also in situations where  $p(x)$ has singularities or is highly peaked in certain regions, Reduced Rejection sampling can be very efficient.  

%In general, PQ-sampling is more efficient than acceptance-rejection method as it rejects less samples. In particular, PQ-sampling is well suited in situations where the function $p(x)$ is changing constantly (see section \ref{section:Comparison} for an elaboration on this point and section \ref{section:Example2} for an example of changing $p(x)$). \\

The next section  describes the Reduced Rejection sampling and proves its validity. Section \ref{section:Comparison} compares Reduced Rejection sampling to other methods (including other generalizations of acceptance-rejection),   highlights advantages of Reduced Rejection sampling in comparison to other methods, and points out some of the challenges in applying Reduced Rejection sampling. In section \ref{section:Example1}, Reduced Rejection sampling is demonstrated on  a simple example. In section \ref{section:Example2},  Reduced Rejection sampling is applied to an example motivated from plasma physics, for which other sampling methods cannot be used efficiently. In section \ref{section:Example3}, we make some comments on how to apply Reduced Rejection in the context of stochastic chemical kinetics. In the appendix, we provide flow charts for the Reduced Rejection algorithm.

\section{Reduced Rejection sampling}\label{section:Reduced Rejection sampling}
Consider a sample space $\Omega$ with Lebesgue measure $\mu$  on $\Omega$, and two functions $q,p:\Omega\rightarrow \mathbb{R}$. Denote
\[I[q]=\int_\Omega q(x)d\mu(x), \qquad I[p]=\int_\Omega p(x)d\mu(s).\]
By sampling from $\Omega$ according to $p(x)$ we mean sampling from $\Omega$ using probability distribution function $p(x)/I[p]$. Partition sample space $\Omega$ into two sets ${\it \mathcal{S}}$ and ${\it \mathcal{L}}$:
\[\mathcal{L}=\{x\in \Omega: p(x)>q(x)\}, \qquad \mathcal{S}=\{x\in \Omega: p(x)\leq q(x)\}.\]
Reduced Rejection sampling is a  method for sampling from $\Omega$ according to $p(x)$ using an auxiliary function $q(x)$. It depends on the following:
\begin{itemize}
\item The values of $I[q]$, $I[p]$ and $\int_{\mathcal{L}}(p(x)-q(x))d\mu(x)$. Note that the last value is needed only for ``Algorithm II'', see subsection \ref{subsection:Algorithm}.
\item A mechanism to sample from $\Omega$ according to $q(x)$.
\item A mechanism to sample from $\mathcal{L}$ according to $p(x)-q(x)$.
\end{itemize}

%Suppose the values of $I[q]$ and $I[p]$ are known. Furthermore, we have a mechanism to sample from $\Omega$ according to $q(x)$ and sample form $\mathcal{L}$ according to $p(x)-q(x)$. In this case, PQ-sampling method provides an efficient way to sample from $\Omega$ according to $p(x)$. 

Whereas the acceptance-rejection method for sampling from $p(x)$ requires a function $q(x)$ that encloses $p(x)$ (i.e., $0\leq p(x) \leq q(x)$ for all $x\in\Omega$),  the
Reduced Rejection sampling algorithm is a generalization of the acceptance-rejection method, that allows $p(x)> q(x)$ for some $x$. 
The Reduced Rejection sampling algorithm is detailed in Section \ref{subsection:Algorithm}, and its validity as a method for sampling  from $\Omega$ according to $p(x)$ is demonstrated in Section \ref{subsection:Proof}.

\subsection{The Reduced Rejection sampling algorithm}\label{subsection:Algorithm}
%We assume existence of a uniform random number generator on interval $(0,1)$. 

The Reduced Rejection sampling method consists of two algorithms (i.e., two different algorithms) depending on the relative values of $I[p]$ and $I[q]$. The outcome of each algorithm is a value $z$ that is an independent sample from $\Omega$ according to $p(x)$.

\textbf{Algorithm I: }$\mathbf{I[p] \geq I[q]}$.

Perform the following steps:
%With probability $(I[p]-I[q])/I[p]$ perform $(i)$; otherwise, perform $(ii)$.

\begin{itemize}
\item[i)] With probability $(I[p]-I[q])/I[p]$, sample $x_0$ from $\mathcal{L}$ according to $p(x)-q(x)$ and accept $z=x_0$.
\item[ii)] Otherwise (with probability $(I[q]/I[p]$), sample $x_0$ from $\Omega$ according to $q(x)$.
  \begin{enumerate}
	\item[a)] If $x_0\in \mathcal{L}$, accept $z=x_0$.
	\item[b)] If $x_0\in \mathcal{S}$,  accept $z=x_0$ with probability $p(x_0)/q(x_0)$. 
  \end{enumerate}
\item[iii)] If $x_0$ was not accepted,  then sample a new value of $x_1$ from $\mathcal{L}$ according to  $p(x)-q(x)$ and accept $z=x_1$.  
\end{itemize}

%\begin{itemize}
%\item[i)] With probability $(I[p]-I[q])/I[p]$ sample from $\mathcal{L}$ according to $p(x)-q(x)$ and accept the result; otherwise,
%\item[ii)] Otherwise, (i.e. with probability $I[q]/I[p]$) sample from $\Omega$ according to $q(x)$; suppose $x_0$ is sampled:
%  \begin{enumerate}
%	\item[a)] If $x_0\in \mathcal{L}$, accept $x_0$.
%	\item[b)] If $x_0\in \mathcal{S}$, with probability $p(x_0)/q(x_0)$ accept $x_0$; otherwise, sample from $\mathcal{L}$ according to function $p(x)-q(x)$ and accept the result.  
%  \end{enumerate}
%\end{itemize}

\textbf{Algorithm II: }$\mathbf{I[p] < I[q]}$.

Perform the following steps until a value $z$ is accepted:

\noindent

\begin{itemize}
\item[i)] Sample $x_0$ from $\Omega$ according to $q(x)$.
%\begin{enumerate}
%\item[a)] If $x_0\in \mathcal{L}$, accept $z=x_0$.
%\item[b)] Sample $x_1$ from $\mathcal{L}$ according to $p(x)-q(x)$ and accept $z=x_1$.
%\item[c)] Return to (i) without accepting a value of $z$.
%\end{enumerate} 
\item[ii)] If $x_0\in \mathcal{L}$, accept $z=x_0$.
\item[iii)] If $x_0\in \mathcal{S}$, accept $z=x_0$ with probability $p_a=p(x_0)/q(x_0)$, 
\item[iv)] If $x_0$ was not accepted, then
\begin{enumerate} 
\item[a)]  With probability  $p_a$ select $x_1$ from $\mathcal{L}$ according to $p(x)-q(x)$ and accept $z=x_1$, in which 
\begin{eqnarray}
p_a &=& \frac{\int_{\mathcal{L}}(p(x)-q(x))d\mu(x)}{\int_{\mathcal{S}}(q(x)-p(x))d\mu(x)}
\nonumber \\
&=&\frac{\int_{\mathcal{L}}(p(x)-q(x))d\mu(x)}{I[q]-I[p]+\int_{\mathcal{L}}(p(x)-q(x))d\mu(x)}.\label{equation:algorithmII} 
\end{eqnarray}
\item[b)] Otherwise (i.e., with probability $1-p_a$), return to (i) without accepting a value of $z$.
\end{enumerate}
\end{itemize}

%\begin{itemize}
%\item[i)] Sample $x_0$ from $\Omega$ according to $q(x)$.
%\item[ii)] If $x_0\in \mathcal{L}$, accept $z=x_0$.
%\item[iii)] Otherwise (i.e., if $x_0\in \mathcal{S}$), select one of  (a) with probability $p_a=p(x_0)/q(x_0)$, or  (b) with probability $p_b$, or  (c) with probability $p_c=1-p_a-p_b$, in which  
%\begin{eqnarray}
%p_b&=& (1-p_a) \frac{\int_{\mathcal{L}}(p(x)-q(x))d\mu(x)}{\int_{\mathcal{S}}(q(x)-p(x))d\mu(x)}
%\nonumber \\
%&=&\frac{q(x_0)-p(x_0)}{q(x_0)}\frac{\int_{\mathcal{L}}(p(x)-q(x))d\mu(x)}{I[q]-I[p]+\int_{\mathcal{L}}(p(x)-q(x))d\mu(x)},\label{equation:algorithmII} 
%\end{eqnarray}
%\begin{enumerate}
%\item[a)] Accept $z=x_0$;
%\item[b)] Sample $x_1$ from $\mathcal{L}$ according to $p(x)-q(x)$ and accept $z=x_1$.
%\item[c)] Return to (i) without accepting a value of $z$.
%\end{enumerate} 
%\end{itemize}

%Note that  performing an action $A$ with probability $r$ (with $0\leq r\leq 1$) is equivalent to choosing a random value $z$ uniformly  from $(0,1)$ and performing $A$ if $z<r$. Also,  {\it accepting} a value $v$  means that $v$ is  the output of the algorithm. 

Figures \ref{figure:AlgorithmI} and \ref{figure:AlgorithmII} in appendix \ref{section:appendix}, illustrate flow charts of Algorithms I and II.

\begin{figure} [!htp]
    \centering
        \includegraphics[width=12cm]{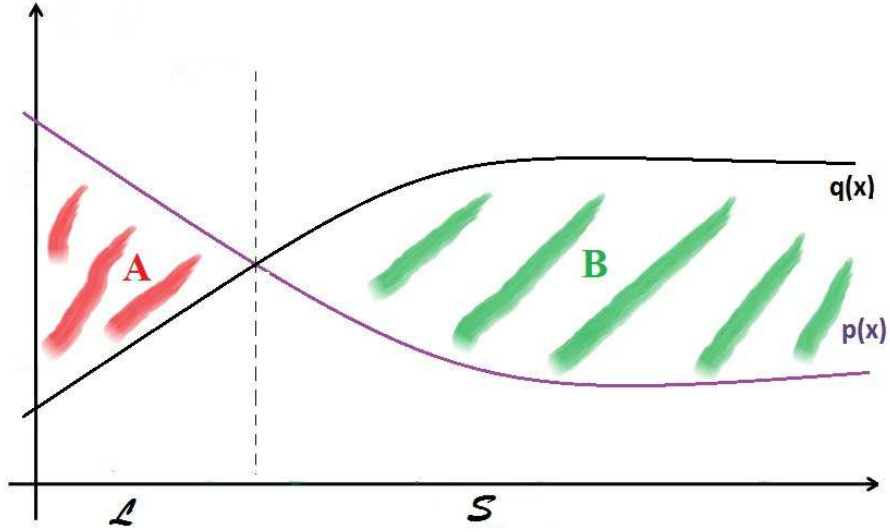} 
    \caption{This figure illustrates the Reduced Rejection method.
Region $A$ is where $q$ does not enclose $p$, and region $B$ is where samples are rejected. Rejected samples from region $B$ can be replaced by samples from region $A$ if $|A|>|B|$; otherwise (if $|A|<|B|$), some of the rejected samples lead to repetition of the algorithm.}
    \label{figure:RRsampling}
\end{figure}

As described in Algorithms I and II, Reduced Rejection samples from $p$ through the following steps:  On $\mathcal{L}$, treat  $p$ as a mixture $p=q+(p-q)$ and sample from $q$ and $p-q$ with the correct probabilities; and on 
$\mathcal{S}$, sample from $p$ by sampling from $q$ and accepting the sample with probability $p/q$. Rejected samples in $\mathcal{S}$ correspond to the region $B$ in figure \ref{figure:RRsampling}, and the region $A$ is where $q$ does not enclose $p$. If $|A|>|B|$ (i.e., Algorithm I) then all of the rejected samples can be replaced by samples from $A$; if $|A|<|B|$ (i.e., Algorithm II) then a portion of the rejected samples can be replaced by samples from $A$, and for the remainder, the algorithm is repeated as in Acceptance-Rejection.

\subsection{Validity of Reduced Rejection Sampling}\label{subsection:Proof}
In this subsection we show the correctness of the Reduced Rejection sampling method. As the method is different for Algorithms I and II, we prove the correctness for each algorithm separately. 

\vspace{8mm}

\textbf{Proof for Algorithm I: }For each $z\in \Omega$, show that the algorithm of Algorithm I returns $z$  with probability $p(z)d\mu(z)/I[p]$.

If $z\in \mathcal{S}$, then  part (ii) must have been selected, $z$ must have been sampled in (ii) and it must have been accepted in case (ii.b).  Therefore, the probability of returning $z$ is 
\begin{eqnarray}
\pr[\hbox{(ii) selected}] \ \pr[z\hbox{ sampled in (ii)}] \ \pr[z \hbox{ accepted in (ii.b)}] 
&=&\frac{I[q]}{I[p]}\times\frac{q(z)d\mu(z)}{I[q]}\times\frac{p(z)}{q(z)}\nonumber \\
&=&\frac{p(z)d\mu(z)}{I[p]}.\label{equation:Ismall}
\end{eqnarray}
Also note that for every $x_0\in \mathcal{S}$, after $x_0$ is selected in (ii.b) with probability $\frac{q(x)}{I[q]}d\mu(x)$, the probability of reaching  (iii) is $\frac{q(x_0)-p(x_0)}{q(x_0)}$. Thus the total probability of reaching  (iii) after selecting (ii) is  
\begin{equation} \pr[\hbox{reaching  (iii) $\vert$  (ii) selected}]= \int_{\mathcal{S}}\frac{q(x)-p(x)}{q(x)}\frac{q(x)}{I[q]}d\mu(x)= \frac{\int_{\mathcal{S}}(q(x)-p(x))d\mu(x)}{I[q]}.\label{equation:IsampleB}\end{equation}

Next suppose that $z\in \mathcal{L}$. The probability that $z$ is returned from  (i)  is
\begin{align}
\pr[z \hbox{ returned from (i)}] =& \pr[\hbox{(i) selected}] \ \pr[z \hbox{ sampled in (i)}] \notag \\
= & \frac{(I[p]-I[q])}{I[p]} \times \frac{(p(z)-q(z))d\mu(z)}{\int_{\mathcal{L}}(p(x)-q(x))d\mu(x)}.
\label{equation:I0}
\end{align}
The probability that $z$ was returned from  (ii.a) is 
\begin{align}
\pr[z \hbox{ returned from (ii.a)}] = & \pr[\hbox{(ii) selected}] \ \pr[z \hbox{ sampled in (ii.a)}]\notag \\
= &\frac{I[q]}{I[p]}\times\frac{q(z)d\mu(z)}{I[q]}\notag \\
= & \frac{q(z)d\mu(z)}{I[p]}.
\label{equation:I1}
\end{align} 
Also, using equation \eqref{equation:IsampleB}, the probability that $z$ was returned from  (iii) is 
\begin{align}
%\begin{array}{ll}
  \pr[z &\hbox{ returned from (iii)}] \notag \\
= & \pr[\hbox{(ii) selected}] \ \pr[\hbox{reaching  (iii) $\vert$  (ii) selected}] \ \pr[z \hbox{ sampled from $\mathcal{L}$ in (iii)}]\nonumber \\ %\end{array}
= & \frac{I[q]}{I[p]}\times\left(\frac{\int_{\mathcal{S}}(q(x)-p(x))d\mu(x)}{I[q]}\right)\times\frac{(p(z)-q(z))d\mu(z)}{\int_{\mathcal{L}}(p(x)-q(x))d\mu(x)} \notag\\  %\end{array}
%\begin{equation}% \begin{array}{ll}
= & \frac{\int_{\mathcal{S}}(q(x)-p(x))d\mu(x)}{I[p]}\frac{(p(z)-q(z))d\mu(z)}{\int_{\mathcal{L}}(p(x)-q(x))d\mu(x)}.
%\end{array} 
\label{equation:I2}
\end{align}
Finally, using equations \eqref{equation:I0}, \eqref{equation:I1}, and \eqref{equation:I2}, the probability of returning $z$ is
\begin{align}
 \pr[z &\hbox{ returned from (i)}]+\pr[z \hbox{ returned from (ii.a)}]+\pr[z \hbox{ returned from (iii)}] \notag\\
=&\frac{(I[p]-I[q])}{I[p]}\frac{(p(z)-q(z))d\mu(z)}{\int_{\mathcal{L}}(p(x)-q(x))d\mu(x)}
+\frac{q(z)d\mu(z)}{I[p]}
+\frac{(p(z)-q(z))d\mu(z)}{\int_{\mathcal{L}}(p(x)-q(x))d\mu(x)}\frac{\int_{\mathcal{S}}(q(x)-p(x))d\mu(x)}{I[p]} \notag\\
=&\frac{(p(z)-q(z))d\mu(z)}{I[p]\int_{\mathcal{L}}(p(x)-q(x))d\mu(x)}\left(I[p]-I[q]+\int_{\mathcal{S}}(q(x)-p(x))d\mu(x)\right)+\frac{q(z)d\mu(z)}{I[p]} \notag\\
=&\frac{(p(z)-q(z))d\mu(z)}{I[p]\int_{\mathcal{L}}(p(x)-q(x))d\mu(x)}\left(\int_{\mathcal{L}}(p(x)-q(x))d\mu(x)\right)+\frac{q(z)d\mu(z)}{I[p]} \notag\\
=&  \frac{p(z)d\mu(z)}{I[p]}.
\label{equation:Ilarge}\end{align}

Hence, by \eqref{equation:Ismall} and \eqref{equation:Ilarge}, whether $z\in \mathcal{S}$ or $z\in \mathcal{L}$, the probability of returning $z$ is equal to $p(z)d\mu(z)/I[p]$. This completes the proof for Algorithm I. \qed

\vspace{8mm}

\textbf{Proof for Algorithm II: }For each $z\in \Omega$, show that  the algorithm in Algorithm II returns $z$ with probability $p(z)d\mu(z)/I[p]$.
The algorithm consists of some number of cycles, each consisting of steps (i)-(iv), until a value $z$ is accepted. We first calculate the probability that $z$ is accepted within one of the cycles.

Suppose that $z\in \mathcal{S}$. Then $z$ must be sampled in (i) and  accepted in (iii). Thus, the probability of returning $z$ in (iii) is 
\begin{align}\pr[z \hbox{ returned from (iii)}]
=&\pr[z \hbox{ sampled in (i)}] \ 
  \pr[z \hbox{ accepted in (iii)}] \notag\\
=&\frac{q(z)d\mu(z)}{I[q]}\times\frac{p(z)}{q(z)}\notag\\
=&\frac{p(z)d\mu(z)}{I[q]}. \label{equation:IIsmall}
\end{align}
Also note that for every $x_0\in \mathcal{S}$,which is chosen with probability $\frac{q(x_0)d\mu(x_0)}{I[q]}$, the probability that it is not accepted in  (iii)  is $\frac{q(x_0)-p(x_0)}{q(x_0)}$. Thus the total probability of not returning an element of $\mathcal{S}$ in  (iii), which is the same as the probability of reaching (iv), is 
\begin{equation} 
\pr[\hbox{reaching (iv)}]
=\int_{\mathcal{S}}\frac{q(x)-p(x)}{q(x)}\frac{q(x)}{I[q]}d\mu(x)
=\int_{\mathcal{S}}\frac{q(x)-p(x)}{I[q]}d\mu(x).  \label{equation:NotReturn}
\end{equation}

Next suppose that $z\in \mathcal{L}$. The probability that $z$ is accepted in (ii) is
\begin{equation} 
\pr[z \hbox{ returned from (ii)}]=\frac{q(z)d\mu(z)}{I[q]}. \label{equation:ai}
\end{equation}
For $z$ to be returned from  (iv.a),  the algorithm must  
reach (iv), then go to (iv.a) and then select $z$ in (iv.a). This has probability
%
%select (b) in (iii) and  sample $z$ in (b). Thus, using equations \eqref{equation:algorithmII} and \eqref{equation:NotReturn}, we have
%
\begin{align}
\pr[z &\hbox{ returned from (iv.a)}] \notag \\
= &     \pr[\hbox{reach (iv)}]  \ 
        \pr[\hbox{go to (iv.a)}]  \ 
        \pr[\hbox{z sampled  in (iv.a)}]  \notag\\ 
= &
\left(\int_{\mathcal{S}}\frac{q(x)-p(x)}{I[q]}d\mu(x)\right)\times
\frac{\int_{\mathcal{L}}(p(x)-q(x))d\mu(x)}{\int_{\mathcal{S}}(q(x)-p(x))d\mu(x)}
\times 
\frac{(p(z)-q(z))d\mu(z)}{\int_{\mathcal{L}}(p(x)-q(x))d\mu(x)} \nonumber \\
= & \frac{(p(z)-q(z))d\mu(z)}{I[q]}. 
\label{equation:bi}
\end{align}
Now using equations \eqref{equation:ai} and \eqref{equation:bi}, the probability of returning $z$ in a cycle is 
\begin{align}
 \pr[z \hbox{ returned}]= & \pr[z \hbox{ returned from (ii)}]+\pr[z \hbox{ returned from (iv.a)}]\notag\\
=& \frac{q(z)d\mu(z)}{I[q]}+\frac{(p(z)-q(z))d\mu(z)}{I[q]}\notag\\
=&\frac{p(z)d\mu(z)}{I[q]}.
\label{equation:IIlarge}
\end{align}
Equations \eqref{equation:IIsmall} and \eqref{equation:IIlarge} imply that, whether $z\in \mathcal{S}$ or $z\in\mathcal{L}$, the probability that $z$ is returned in a cycle is $p(z)d\mu(z)/I[q]$. Integrating over all the samples in $\Omega$, we deduce that the probability that a sample is returned in a cycle is $I[p]/I[q]$. Consequently, probability that no sample point is returned in a cycle is $1-I[p]/I[q]$. Because the cycle is repeated until a sample point is returned, we conclude that the probability that the algorithm returns $z$ is equal to 
\[\sum_{k=1}^{\infty}(1-\frac{I[p]}{I[q]})^{k-1}\frac{p(z)d\mu(z)}{I[q]}=\frac{p(z)d\mu(z)}{I[p]}.\]   
This completes the proof for Algorithm II.\qed

Note that the efficiency of Algorithm II is nominally the same as acceptance rejection, i.e. the probability of a rejection is $1- I[p]/I[q]$. Actually it can be significantly better because $I[q]$ can be smaller, since $q<p$ is allowed. 
Also, note that if $I[p]=I[q]$, then Algorithms I and II are the same. 

\section{Comparison of Reduced Rejection  and Other Sampling Methods}\label{section:Comparison}
%In sampling a random variable, the choice of the algorithm depends on different parameters and in certain situations some methods perform better than other methods. As mentioned earlier in this paper, we have focused on methods that provide independent samples. 
One of the important features of Reduced Rejection sampling is that it requires no preprocessing time. This is  particularly  useful for {\it dynamic simulation}; i.e., simulation in which the probability distribution function $p(x)$ may change after each sample (see section \ref{section:Example2} for an example from plasma physics). For dynamic simulation, fast discrete sampling methods such as Marsaglia's table method or the Alias method, are not suitable as they require preprocessing time after each change in  $p(x)$. Although, the acceptance-rejection method requires no preprocessing time and can be used for dynamic simulation, it may require changes in $q(x)$ if $p(x)$ changes, which is usually not difficult, and it becomes very inefficient when the ratio of the area under function $p(x)$ to the area under proposal function $q(x)$ is small. Moreover, adaptive rejection sampling  is not  efficient,  because the process of adapting  $q(x)$ to $p(x)$ starts over whenever $p(x)$ changes.

The Reduced Rejection sampling method can be thought as an extension of the acceptance-rejection method. In particular when the proposal function $q(x)$ encloses $p(x)$ (i.e., 
$q(x)\geq p(x) \hbox{ for all } x\in\Omega$ so that $\mathcal{L}=\emptyset$) the Reduced Rejection sampling method reduces to acceptance-rejection method. The advantage of Reduced Rejection sampling over acceptance-rejection method is that the proposal function $q(x)$ does not  need to enclose function $p(x)$; i.e., it allows $q(x)<p(x)$ for some $x$. This is very useful in dynamic simulation as it can accommodate changes in $p(x)$ without requiring changes in  $q(x)$. Moreover, Reduced Rejection sampling may result in less unwanted samples  than acceptance-rejection does, especially if  $p(x)$ has singularities or is highly peaked. 
%Moreover, in general PQ-sampling rejects less samples than acceptance-rejection does. Indeed the probability of rejecting a sample in Situation II of PQ-sampling is $(I[q]-I[p])/I[q]$, which is the same formula as for the acceptance-rejection. However, because functions $q(x)$ and $p(x)$ do not need to satisfy condition \eqref{equation:criteria} in PQ-sampling, in practice, the number of times we need to sample from proposal function $q(x)$ to generate specific number of independent samples according to $p(x)$ is less than it is in  
%the probability of rejecting a sample in PQ-sampling is less than it is in acceptance-rejection. \\

There are several challenges in implementing the Reduced Rejection sampling method. The main challenge is the need to sample from set $\Omega$ according to $q(x)$ and from set $\mathcal{L}$ according to $p(x)-q(x)$, which can be performed by various sampling methods. 
%Depending on the problem, one uses different sampling methods, for example acceptance-rejection, to complete this task. Consider a function $p(x)$ that is highly peaked in certain regions of the domain. We choose our proposal function $q(x)$ in such a way that set $\mathcal{L}$ becomes the regions at which $p(x)$ is highly peaked. Sampling from $\mathcal{L}$ using acceptance-rejection is far more efficient than sampling from the whole domain using acceptance-rejection. Therefore, PQ-sampling boosts the sampling process.
Another challenge in using Reduced Rejection sampling is the need to know the values of $I[q]$, $I[p]$ and $\int_{\mathcal{L}}(p(x)-q(x))d\mu(x)$ (but note that the last value is only for Algorithm II). In many situations, these values are readily available or can be calculated during the simulation. 

\section{Example 1: Reduced Rejection Sampling for a Random Variable with Singular Density}\label{section:Example1}
In this section, Reduced Rejection sampling method is applied to a simple problem. Let $\Omega=(0,1)$ and  sample according to  
\begin{equation}p(x)=\frac{1}{\sqrt{x}}+\frac{1}{\sqrt[5]{1-x}}\label{equation:Ex1}\end{equation}
which has singularities at 0 and 1. Using inverse transform sampling, it is  easy to sample according to $1/\sqrt{x}$ or $1/\sqrt[5]{1-x}$, but inverse transform cannot be easily applied to  \eqref{equation:Ex1} as it requires finding the root of a eighth degree polynomial. We apply Reduced Rejection sampling to this problem by setting  $q(x)=1/\sqrt{x}$. Observe that $\mathcal{L}=\Omega=(0,1)$. As mentioned earlier, inverse transform sampling is easily used to sample  according to $q(x)$ and  according to $p(x)-q(x)$. The Reduced Rejection sampling is very fast and yields no unwanted sample points. 
%Although  other methods can also be used to generate samples efficiently;  none of them seems to be as fast and efficient as Reduced Rejection sampling for this problem.  
This example is equivalent to sampling from a mixture and can be extended to  sampling from a probability density $p(x)$ that is a sum $p = p_1 + p_2 + \ldots + p_n$, if there is a method for sampling from each $p_k$ separately and the integrals $I[p_k]$ are all known.

\section{Example 2: Reduced Rejection Sampling for a Stochastic Process with Fluctuating and Singular Rates}\label{section:Example2}
In this section, we apply Reduced Rejection sampling to an idealized problem motivated by plasma physics. As discussed in subsection \ref{subsection:Numerical result}, the unique features of this problem makes other sampling methods  inefficient to use.

\subsection{Statement of the Stochastic Process and the Simulation Algorithm }\label{subsection:Statement}

The example presented here is a simplified version of simulation for recombination by impact of two electrons with an ion, in which one of the electrons is absorbed into the atom and the other electron is scattered. For incident electron energies $E_1$ and $E_2$, the recombination rate is proportional to $(E_1 E_2)^{-1/2}$ \cite{BIB:Oxenius,BIB:Zeldovich}, which can become singular if  electrons of low energy are present. This is an obstacle to kinetic simulation of  recombination by electron impact in a plasma.

Our goal is to simulate the evolution of the following system:  Consider $N$ particles labeled $1,\ldots,N$. To each particle $i$ we associate a number $x_i\in (0,1)$, called the state of particle $i$ (and corresponding to electron energy in the recombination problem). Occasionally, where it does not cause confusion, we use $x_i$ to refer to particle $i$. We refer to the set $\Gamma=\{x_1,\ldots,x_n\}$ as the configuration of the system. For every pair of states $x_i$ and $x_j$,   $T_{i,j}$ is a random variable with an exponential distribution with parameter $(x_ix_j)^\alpha$, in which $\alpha$ is a fixed constant between 0 and 1.  $T_{i,j}$ is the  time for interaction between particle $i$ and $j$ which randomly occurs with rate $1/(x_ix_j)^\alpha$. After an interaction occurs, say for the pair $\{k,l\}$,  the values of states $x_k$ and $x_l$ are replaced by new values $x'_k$ and $x'_l$; consequently, the distribution of $T_{i,j}$ changes if either of $i$ and $j$ is equal to $k$ or $l$. 

 We will consider a simple updating mechanisms for the states after each interaction. In the simulations presented below, the updated values of $x'_k$ and $x'_l$ are chosen independently and uniformly at random from $(0,1)$, without dependence on $x_k$ and $x_l$. This choice has been made for simplicity and because the stationary distribution can be calculated for this choice 
(see section \ref{subsection:Theoretical result}), but we expect that Reduced Rejection sampling would work equally well for more complex interaction rules.
%in the second they are chosen randomly subject to the condition that ${x'_k}^2+{x'_l}^2=x_k^2+x_l^2$.
Indeed the  algorithm \ref{algorithm:process} described below and the more detailed algorithm presented in section \ref{subsection:Use of Reduced Rejection } do not depend on the interaction rules.

First we make some notation and observations. Set $s_i=1/x_i^\alpha$ and $s=\sum_{i}s_i$. Let $T(\lambda)$ denotes an exponential random variable with parameter $\lambda$ (with rate $1/\lambda$); then $T(\lambda)=\mu T(\mu\lambda)$ for any scalar $\mu$. We will use
\[\frac{s_i}{s}\frac{s_j}{s}T(1/s^2)=T(1/(s_is_j))\]
in  the following algorithm, which is a variant of the Kinetic Monte Carlo (KMC) algorithm (also known as the residence-time algorithm or the n-fold way or the Bortz-Kalos-Lebowitz (BKL) algorithm \cite{BIB:Bortz1975}), that simulates the system described above. This algorithm chooses interactions, by choosing two particles separately out of the  $N$  number of particles, rather than choosing a pair of particles out of the $N^2$ number of pairs.
\begin{algorithm}\label{algorithm:process}

\begin{enumerate}
\item Start from $t=0.$
\item Choose time $\Delta t$ by sampling from an exponential distribution with rate $s^2$.
\item Choose index $k$ with probability $s_k/s$.
\item Choose index $l$ with probability $s_l/s$.
\item At time $t+\Delta t$ interaction between particles $k$ and $l$ occurs.
\item Update states $x_k$ and $x_l$ according to the updating mechanism and the update value of $s$. 
\item Set $t=t+\Delta t$ and start over from 2. 
\end{enumerate}
\end{algorithm}

We use Reduced Rejection sampling in subsection \ref{subsection:Numerical result} to perform steps 3 and 4 in the above algorithm. We also explain why other methods of sampling would be inefficient in this circumstances. To verify that our simulation is working properly, we perform the following test.

Let $g(x_1,\ldots,x_N)$ be a real-valued function on configuration space, with expectation $E[g]$ of  $g$ over configurations of the system. For a simple  interaction rule  and some functions $g$ we can find the value of $E[g]$ analytically, as shown in subsection \ref{subsection:Theoretical result}. Consequently, the difference between the numerical and analytic results provides a measure of the accuracy of the simulation as discussed at the end of subsection \ref{subsection:Numerical result}.

\subsection{Theoretical results}\label{subsection:Theoretical result}
Think of the system's evolution  as a random walk over the configurations of the system. Suppose the updating process is that if states $x_k$ and $x_l$ interact, then states $x'_k$ and $x'_l$ are chosen, independently, uniformly at random from $(0,1)$. In this section, we find the stationary distribution for this random walk and the value of $E[g]$ for two functions $g$.

 Let $P(\Gamma'|\Gamma)$ denote the probability of going from configuration $\Gamma$ to configuration $\Gamma'$, and $P(x'_k,x'_l|x_k,x_l)$ as the probability of going from values $x_k,x_l$ to values $x'_k,x'_l$ in an interaction.  These satisfy 
\begin{equation} P(x'_k,x'_l|x_k,x_l)=P(x_k,x_l|x'_k,x'_l),\label{equation:probRelation}\end{equation}
that is, the probability of getting $x'_k$ and $x'_l$ after $x_k$ and $x_l$ interact, is the same as probability of getting $x_k$ and $x_l$ after $x'_k$ and $x'_l$ interact. Furthermore, for this updating mechanism, the random walk is completely mixing; that is, it can go from any configuration to any other configuration (If for example the updating mechanism had additional constraints, such as ${x'_k}+{x'_l}=x_k+x_l$, then we would not have a mixing random walk since we could reach only those configurations that have the same sum of states as the starting configuration).

For every configuration $\Gamma=\{x_1\ldots,x_N\}$, set
\[   \pi_\Gamma:=\frac{(x_1\cdots x_N)^\alpha}{Z}\left(\sum_{i,j}\frac{1}{(x_ix_j)^\alpha}\right),\] 
where $Z$ is the normalizing constant so that $\sum_\Gamma \pi_\Gamma=1$.

Suppose the current configuration of the system is $\Gamma=\{x_1,\ldots,x_N\}$. According to steps 3 and 4 in algorithm \ref{algorithm:process}, the probability of interaction occurring between states $x_k$ and $x_l$ is proportional to $s_{k}s_{l}=1/(x_{k}x_{l})^\alpha$. Let $\Gamma'=\{x'_k,x'_l\}\cup\Gamma\setminus\{x_k,x_l\}$. Then
\[P(\Gamma'|\Gamma)=\frac{1/(x_kx_l)^\alpha}{\sum_{i,j}1/(x_ix_j)^\alpha}P(x'_k,x'_l|x_k,x_l).\]
For the ease of explanation, relabel the states of $\Gamma'$ so that $\Gamma'=\{x'_1,\ldots,x'_N\}$ where $x'_i=x_i$ for $i\neq k,l$. Similarly,
\[ \pi_{\Gamma'}=\frac{(x'_1\cdots x'_N)^\alpha}{Z}\left(\sum_{i,j}\frac{1}{(x'_ix'_j)^\alpha}\right), \quad\hbox{ and }\quad P(\Gamma|\Gamma')=\frac{1/(x'_kx'_l)^\alpha}{\sum_{i,j}1/(x'_ix'_j)^\alpha}P(x_k,x_l|x'_k,x'_l). \]
Because of \eqref{equation:probRelation} and since $x'_i=x_i$, for $i\neq k,l$, it is straightforward to verify the detailed balance equation
\[\pi_\Gamma P(\Gamma'|\Gamma)=\pi_{\Gamma'} P(\Gamma|\Gamma').\]
Therefore, $\pi_\Gamma$ is the (unique) stationary distribution of the random walk. Since the updating mechanism is completely mixing, the normalizing constant $Z$ for distribution $\pi_\Gamma$ is the integral
\[ Z=\int_{(0,1)^N}(x_1\cdots x_N)^\alpha\left(\sum_{i,j}\frac{1}{(x_ix_j)^\alpha}\right)dx_1\cdots dx_N=\frac{\binom{N}{2}}{(\alpha+1)^{N-2}}.\]
Hence for any function $g(x_1,\ldots,x_N)$, 
\[E[g]=\frac{(\alpha+1)^{N-2}}{\binom{N}{2}}\int_{(0,1)^N}g(x_1,\ldots,x_N)(x_1\cdots x_N)^\alpha\left(\sum_{i,j}\frac{1}{(x_ix_j)^\alpha}\right)dx_1\cdots dx_N.\]
Some tedious algebra leads to the following proposition:
\begin{proposition}\label{proposition:g}
Using the above notations and assumptions:
\begin{itemize}
\item[a)] $E[g]=\frac{\alpha+1}{\alpha+2}(N-2)+1$ when $g(x_1,\ldots,x_N)=x_1+\cdots+x_N$.
\item[b)] $E[g]=\frac{\alpha+1}{\alpha+3}(N-2)+\frac{2}{3}$ when $g(x_1,\ldots,x_N)=x^2_1+\cdots+x^2_N$.
\end{itemize}
\end{proposition}

\subsection{Simulation issues}\label{subsection:Simulation issues}
In this section we  make some remarks about the challenges involved in simulating this system.

The main challenge of sampling in this  dynamic simulation is that the  $s_i$'s are changing after each interaction. Consequently, the sampling method should require small or zero preprocessing time. For this reason, discrete sampling methods such as Marsaglia's table method or the Alias method are not very efficient for this problem.

Next  consider using acceptance-rejection method based on uniform sampling from $1$ to $n$ for the proposal distribution (i.e., $q$ constant). As mentioned earlier, the changing distribution property of the problem is not very detrimental for acceptance-rejection. 
On the other hand, the singularity in the rates at $x_i=0$ can lead to a large constant for $q$, for which there will be many rejected samples, so that the method is inefficient.
%However, note that because of the way $s_i$'s are defined, intuitively, this is similar to using acceptance-rejection with uniform proposal distribution to sample according to $1/x^\alpha$ on $(0,1)$. It is clear that there would be far too many rejected samples and the method would be inefficient. 
Moreover, there seems to be no other clear choice for the proposal distribution $q$ other than a constant. Note that the sampling is from a discrete set of probabilities $s_i/s$ with little control over their values; for example the $s_i$'s are not  monotonically ordered. This is quite different from sampling a single random variable from the density $p(x)=x^{-\alpha}$.

%to be used with acceptance-rejection method\footnote{It is also useful to keep in mind that $s_i$'s are not necessarily ordered monotonically in the memory and as they are changing after each interaction, ordering them would require $O(N)$ operations after each interaction, which is costly.}.

\subsection{Use of the Reduced Rejection algorithm}\label{subsection:Use of Reduced Rejection }

In this section we explain how to use Reduced Rejection Sampling to perform steps 3 and 4 in algorithm \ref{algorithm:process}. Reduced Rejection sampling can be readily used in this dynamic simulation. Even though the values of the $s_i$'s change after each interaction, they do not change drastically; in each interaction at most two of the $s_i$'s change. Starting at time $0$, we set  $q_i=p_i=s_i$. After each interaction, we update the values of $p_i$'s to  $p_i=s_i$, but do not change the values of $q_i$'s. Note that we can easily update the value of $I[p]$ after each interaction and keep track of set $\mathcal{L}=\{i:p_i>q_i\}$ by comparing the updated values of $p_i$'s to their corresponding values of $q_i$'s. Moreover, the size of set $\mathcal{L}$ changes by at most 2 after each interaction (but it can also decrease after some interactions). 

We use Marsaglia's table method to sample according to $q_i$'s. Since we do not update $q_i$'s after each interaction, the preprocessing time in Marsaglia's table method is only required for the first sampling and not for the subsequent samplings. To sample from set $\mathcal{L}$ according to $p_i-q_i$, we use acceptance-rejection with uniform distribution for the proposal distribution. As long as the size of set $\mathcal{L}$ is not too big, the sampling from $\mathcal{L}$ is not  very time consuming. To prevent $\mathcal{L}$ from getting too large, we reset the values of $q_i$'s to  $q_i=p_i=s_i$, which sets  $\mathcal{L}$ to be empty, whenever the size of $\mathcal{L}$ exceeds a predetermined number $M$.

The size of $M$ is important for the performance of the algorithm. If  $M$ is too small, then there are many updates of the $q_i$'s, each of which requires preprocessing time for Marsaglia's table method. On the other hand, if  $M$ is too big, then  $\mathcal{L}$ is large and  costly to sample from by acceptance-rejection. Our computational experience shows that setting $M$ equal to a multiple of $\sqrt{N}$ is a good choice.
It might be better for the reinitialization criterion to be based on the efficiency of the sampling from  $\mathcal{L}$ (i.e., the fraction of rejected samples when using acceptance rejection on  $\mathcal{L}$), rather than the size of  $\mathcal{L}$.  

\subsection{Numerical result}\label{subsection:Numerical result}

\begin{figure} [!htp]
    \centering
        \includegraphics[width=12cm]{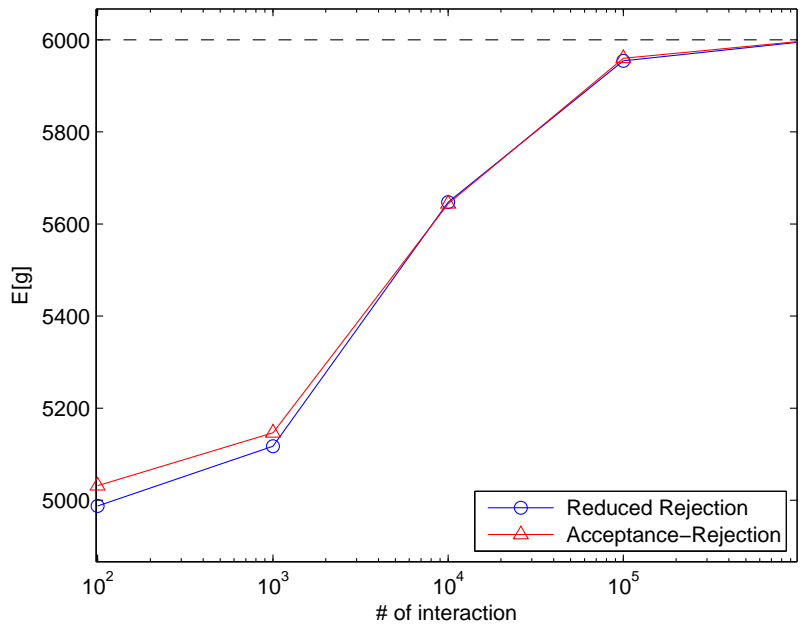} 
    \caption{Theoretical (dashed line) and estimated values (solid lines) of $E[g]$ using different number of interactions. Here $g(x_1,\ldots,x_N)=x_1+\cdots+x_N$, $N=10^4$, and $\alpha=0.5$. Also $M=4000$ for the Reduced Rejection sampling. The theoretical value of $E[g]$ is 5999.8. The estimated value of $E[g]$ after $10^6$ interactions using Reduced Rejection sampling and acceptance-rejection methods were, respectively, 5994.59 and 5996.35. The reported result is the average of 5 independent runs.}
    \label{figure:linear10000result}
\end{figure}

We simulated the evolution of the system under the conditions outlined in subsection \ref{subsection:Theoretical result} with 
 $N=10^4$, and $\alpha=0.5$. We start with a random configuration at time $t=0$.  The simulation is based the Reduced Rejection sampling method, using  Marsaglia's table method and the acceptance-rejection method as described above.  After each interaction, we evaluate function $g(x_1,\ldots,x_N)=x_1+\cdots+x_N$ and take the average to get an estimate for $E[g]$. Each result is produced by taking an average of five independent runs. Figure \ref{figure:linear10000result} compares the results for $E[g]$ from Reduced Rejection sampling with those from the acceptance-rejection method. The results of figure \ref{figure:linear10000result} show excellent agreement between the values of $E[g]$ as a function of the number of interactions from the two methods, which provides a validity check for  Reduced Rejection sampling.

\begin{figure} [!htp]
    \centering
        \includegraphics[width=12cm]{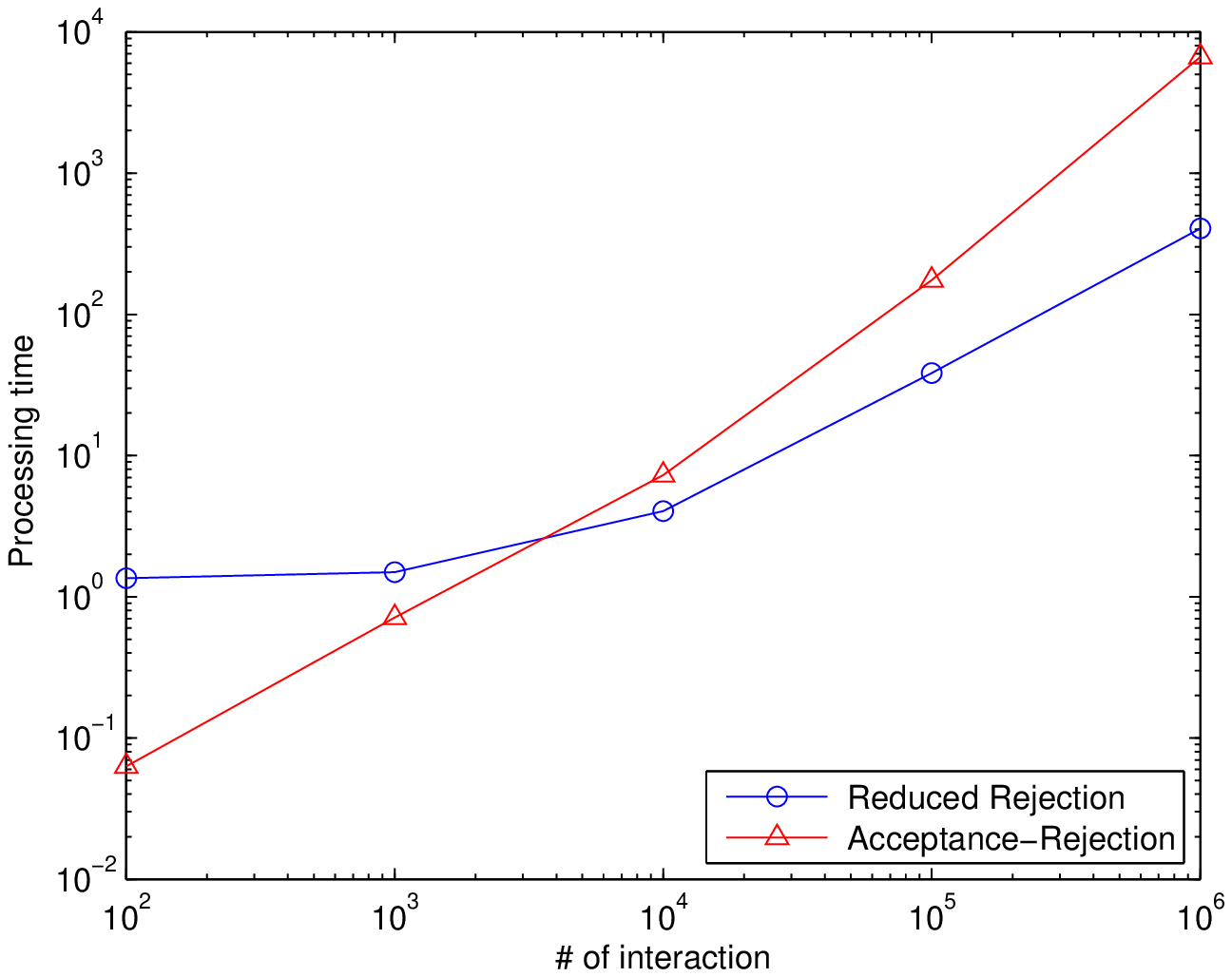} 
    \caption{Loglog plot of the processing time for acceptance-rejection and Reduced Rejection sampling. Here $g(x_1,\ldots,x_N)=x_1+\cdots+x_N$, $N=10^4$ and $\alpha=0.5$. Also $M=4000$ for the Reduced Rejection sampling. The reported processing time is the average of 5 independent runs.}
    \label{figure:linear10000loglog}
\end{figure}

The advantage of Reduced Rejection sampling is demonstrated in figure \ref{figure:linear10000loglog} which shows a log-log plot of the processing time, as a function of the number of interactions, for Reduced Rejection sampling and acceptance-rejection. The results show that Reduced Rejection sampling is much faster than the acceptance-rejection method. In fact, for $n$ interactions, the computational time scales like $O(n)$ for Reduced Rejection  sampling, and like $O(n^{3/2})$ for acceptance-rejection, in the range $10^4 \leq n \leq 10^6$. For small values of $n$, the initial pre-processing step of Marsaglia's table method dominates the computational time. For $n>10^4$, however, the pre-processing time (including the multiple pre-processing steps due to reinitialization) is not a significant part of the computational time. The average  number of reinitializations for Reduced Rejection sampling is $(0,0,0,3.7,53.1)$ for $n=(10^2,10^3,10^4,10^5,10^6)$, respectively. Another interesting advantage of the Reduced Rejection sampling  is that the variance of the processing time for independent runs is much smaller in the Reduced Rejection sampling than it is in the acceptance-rejection method. 

%\begin{table}[!htp] \centering
%\begin{tabular}{|c|c|c|}
%\hline
%\hbox{number of interactions} & E[g]\hbox{ (Reduced Rejection sampling)} & E[g]\hbox{ (acceptance-rejection)} \\ \hline
% $10^2$		&  4987.652463 	&	5031.518436	\\ \hline
% $10^3$		& 5117.149268; 	&	5146.362257 	\\ \hline
% $10^4$		& 5647.174861; 	&	5643.942558	\\ \hline
% $10^5$		& 5954.401000 	&	 5959.757596 	\\ \hline
% $10^6$		& 5994.590683 	&	5996.349960	\\ \hline 
%\end{tabular}
%\caption{The estimated result for $E[g]$, where $g(x_1,\ldots,x_N)=x_1+\cdots+x_N$, by taking the average of value of $g$ after each interaction. The reported result is the average of 5 independent runs. The theoretical value of $E[g]$ is 5999.8 for $N=10000$.}
%\label{table:linearG}\end{table}

%%%%%%%%%%%%%%%%%%%%%%%%%%%%%%
% put back in the following  figures
%%%%%%%%%%%%%%%%%%%%%%%%%%%%%%%%%%%%

\section{Example 3: Stochastic Simulation of Chemical Kinetics}\label{section:Example3}
In this section we describe how we can use Reduced Rejection in the context of stochastic chemical kinetics. Stochastic simulation in chemical kinetics is a Monte Carlo procedure to numerically simulate the time evolution of a well-stirred chemically reacting system. The first Stochastic Simulation Algorithm, called the Direct Method, was presented in \cite{BIB:Gillespie76}. The Direct Method is computationally expensive and there have been many adaptations of this algorithm to achieve greater speed in simulation. The first-reaction method, also in \cite{BIB:Gillespie76}, is an equivalent formulation of the Direct Method. The next-reaction method \cite{BIB:Gibson} is an improvement over the first-reaction method, using a binary-tree structure to store the reaction times. The Modified Direct Method \cite{BIB:Cao} and Sorting Direct Method \cite{BIB:McCollum} speed up the Direct Method by indexing the reactions in such a way that reactions with larger propensity function tend to have a lower index value. Recently, some new Stochastic Simulation Algorithms, called partial-propensity methods, were introduced that work only for elementary chemical reactions (i.e. reactions with at most two different reactants)  (see \cite{BIB:Ramaswamy09,BIB:Ramaswamy10,BIB:Ramaswamy11}). Nevertheless, note that it is possible to decompose any non-elementary reaction into combination of elementary reactions. There are also approximate Stochastic Simulation Algorithms, such as tau-leaping and slow-scale, that provide better computational efficiency in exchange for sacrificing some of the exactness in the Direct Method (see \cite{BIB:Gillespie07} and the references therein for more details).

Next we give a brief review of stochastic simulation in chemical kinetics. An excellent reference with more detailed explanation is \cite{BIB:Gillespie07}. Using the same notation and terminology as in \cite{BIB:Gillespie07}, consider a well-stirred system of molecules of $N$ chemical species $\{S_1,\ldots, S_N\}$, which interact through $M$ chemical reactions $\{R_1,\ldots,R_M\}$. Let $X_i(t)$ denote the number of molecules of species $S_i$ in the system at time $t$. The goal is to estimate the state vector $\X(t)\equiv (X_1(t),\ldots,X_N(t))$ given the system is initially in state $\X(0)=\x_0$.

Similar to section \ref{section:Example2}, when the system is in state $\x$, the time for reaction $R_j$ to occur is given by an exponential distribution whose rate is the propensity function $a_j(\x)$. When reaction $R_j$ occurs, the state of the system changes from $\x$ to $\x+(v_{1j},\ldots, v_{Nj})$, where $v_{ij}$ is the change in the number of $S_i$ molecules when one reaction $R_j$ occurs.

Estimating the propensity functions in general is not an easy task. As noted, the value of the propensity functions depend on the state of the system. For example, if $R_i$ and $R_j$ are, respectively, the unimolecular reaction $S_1\rightarrow product(s)$ and bimolecular reaction $S_1+S_2\rightarrow product(s)$, then $a_i(\x)=c_ix_1$ and $a_j(\x)=c_jx_1x_2$ for some constants $c_i$ and $c_j$. Therefore, the propensity functions of the reactions are changing throughout the simulation. Moreover, if for some chemical species the magnitude of their population differ drastically from others, we expect the value of propensity functions to be very non-uniform.

For every state $\x$, define 
\[  a(\x)=\sum_{i=1}^M a_i(\x).\]
To simulate the chemical kinetics of the system the following algorithm is used, which resembles algorithm \ref{algorithm:process} in section \ref{section:Example2}.
\begin{algorithm}\label{algorithm:chemicalProcess}

\begin{enumerate}
\item Start from time $t=0$ and state $\x=\x_0$.
\item Choose time $\Delta t$ by sampling from an exponential distribution with rate $a(\x)$.
\item Choose index $k$ with probability $a_k(\x)/a(\x)$.
\item At time $t+\Delta t$ reaction $R_k$ occurs.
\item Update time $t=t+\Delta t$, state $\x=\x+(v_{1k},\ldots,v_{Nk})$ and start over from 2.
\end{enumerate}
\end{algorithm}

 In the original Direct Method \cite{BIB:Gillespie76} step 3 in the above algorithm \ref{algorithm:chemicalProcess} is performed by choosing number $r$ uniformly at random in the unit interval and setting
\begin{equation} k=\hbox{ the smallest integer satisfying } \sum_{i=1}^k a_i(\x)>r a(\x).\label{equation:partialSum}\end{equation}
However, when we have many reactions with a wide range of propensity function values presented in the system, a scenario that is very common in biological models, the above procedure of using partial sums becomes computationally expensive. As noted earlier, some methods, such as the Modified Direct Method \cite{BIB:Cao} and Sorting Direct Method \cite{BIB:McCollum}, index the reactions in a smart way so that they can save on the average number of terms summed in equation \eqref{equation:partialSum} and consequently achieve computational efficiency.

We propose a different approach to performing step 3 in algorithm \ref{algorithm:chemicalProcess} using the acceptance-rejection or Reduced Rejection method. The approach is very similar to what was done in section \ref{section:Example2}. To be specific, we can use acceptance-rejection for step 3 in the following way: let 
\[ \bar{a}(\x)=\max_{i} a_i(\x). \]
Until an index is accepted, select index $k$ uniformly at random from $\{1,\ldots,N\}$ and accept it with probability $a_k(\x)/\bar{a}(\x)$; otherwise, discard $k$ and repeat. When an index is accepted step 3 in algorithm \ref{algorithm:chemicalProcess} is completed. Typically for chemical reactions $R_j$, most of $v_{ij}$'s are zero; therefore, we can efficiently update the value of $\bar{a}(\x)$ at each iteration of algorithm \ref{algorithm:chemicalProcess}.

However, as  in section \ref{section:Example2}, if the values of $a_i(\x)$'s are very non-uniform (for example, when the population of some chemical species differ drastically from  that of other species in the system) the acceptance-rejection method becomes inefficient due to rejection of many samples. In these circumstances, the Reduced Rejection algorithm can be readily used in a very similar way as it was used in section \ref{section:Example2}. We expect that the use of the Reduced Rejection algorithm in these circumstances would greatly improve the computational efficiency of the exact Stochastic Simulation Algorithms.

\section{Conclusions and Future Directions}\label{section:Conclusion}
In this paper we introduce a new Reduced Rejection sampling method that can be used to generate independent samples for a discrete or continuous random variable. The strength of this algorithm is most evident for applications in which acceptance-rejection method is inefficient; namely, the probability distribution of the random variable is highly peaked in certain regions or has singularities. It is also useful when the probabilities are fluctuating, so that  discrete methods that requiring preprocessing are inefficient. In particular, the Reduced Rejection sampling method is expected to perform well on kinetic simulation of electron-impact recombination in a plasma, which is difficult to simulate by other methods.

The preliminary examples  in this paper are meant to illustrate these advantages of the Reduced Rejection sampling method. They provide evidence of improvement in computation time using the Reduced Rejection sampling versus acceptance-rejection method. These examples also provide some insights on implementation of the method.
%One limitation in using PQ-sampling is that it requires sampling according to $p(x)-q(x)$ on the set $\mathcal{L}$. Depending on the problem, we can use different methods for this part of the algorithm (for example for the example in section \ref{section:Example1} we used inverse transform sampling and for the example in section \ref{section:Example2} we used acceptance-rejection method).  

One possible direction for future research is  the nested use of  Reduced Rejection sampling methods. For the most difficult step - sampling from   $\mathcal{L}$ according to $p(x)-q(x)$ - we propose to apply the Reduced Rejection sampling method again using  a new proposal function. In essence, this would use one Reduced Rejection sampling method inside  another Reduced Rejection sampling method.  

\appendix

\section{Flow Charts of the Reduced Rejection Algorithm}\label{section:appendix}

\begin{figure} [!htp]
    \centering
        \includegraphics[width=16cm]{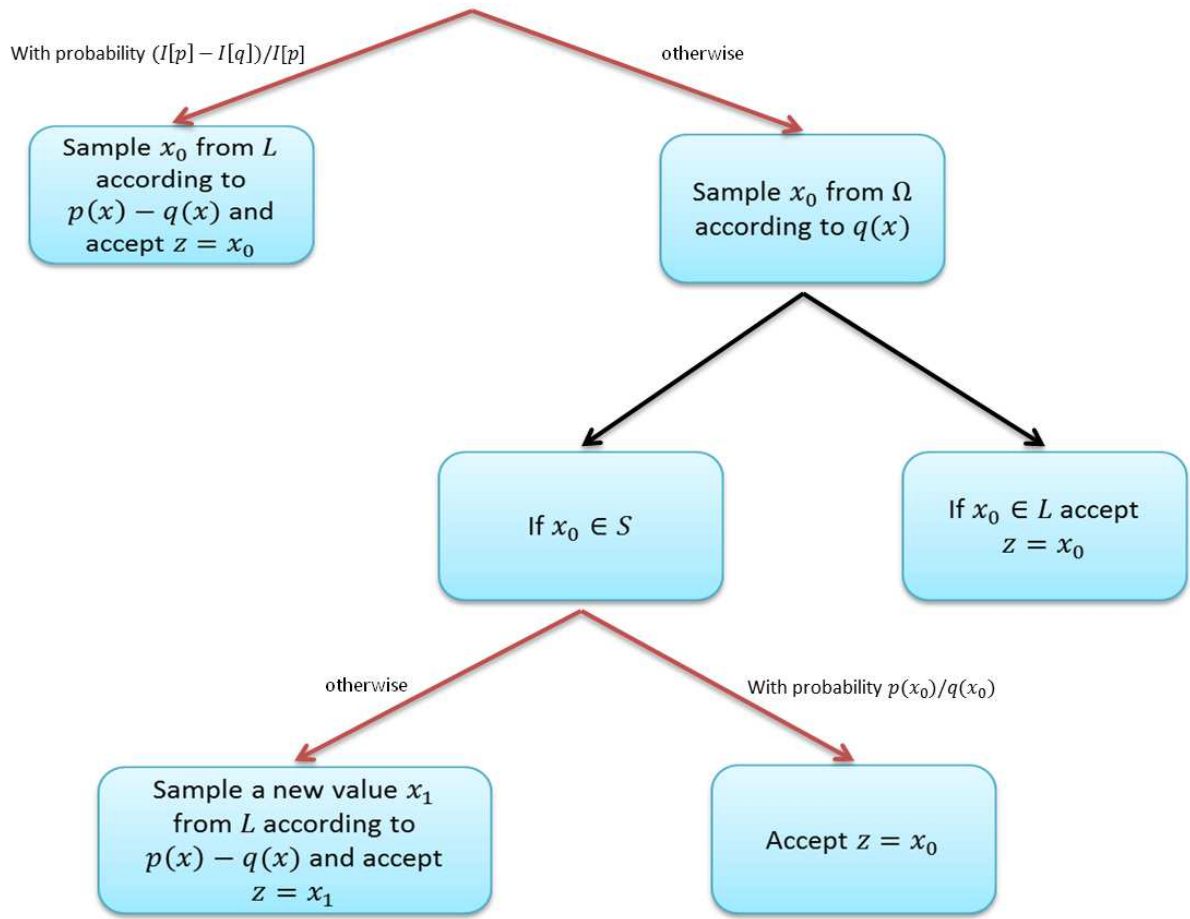} 
    \caption{Flow chart of Algorithm I of the Reduced Rejection sampling method attributing to the case $I[p] \geq I[q]$.}
    \label{figure:AlgorithmI}
\end{figure}

\begin{figure} [!htp]
    \centering
        \includegraphics[width=16cm]{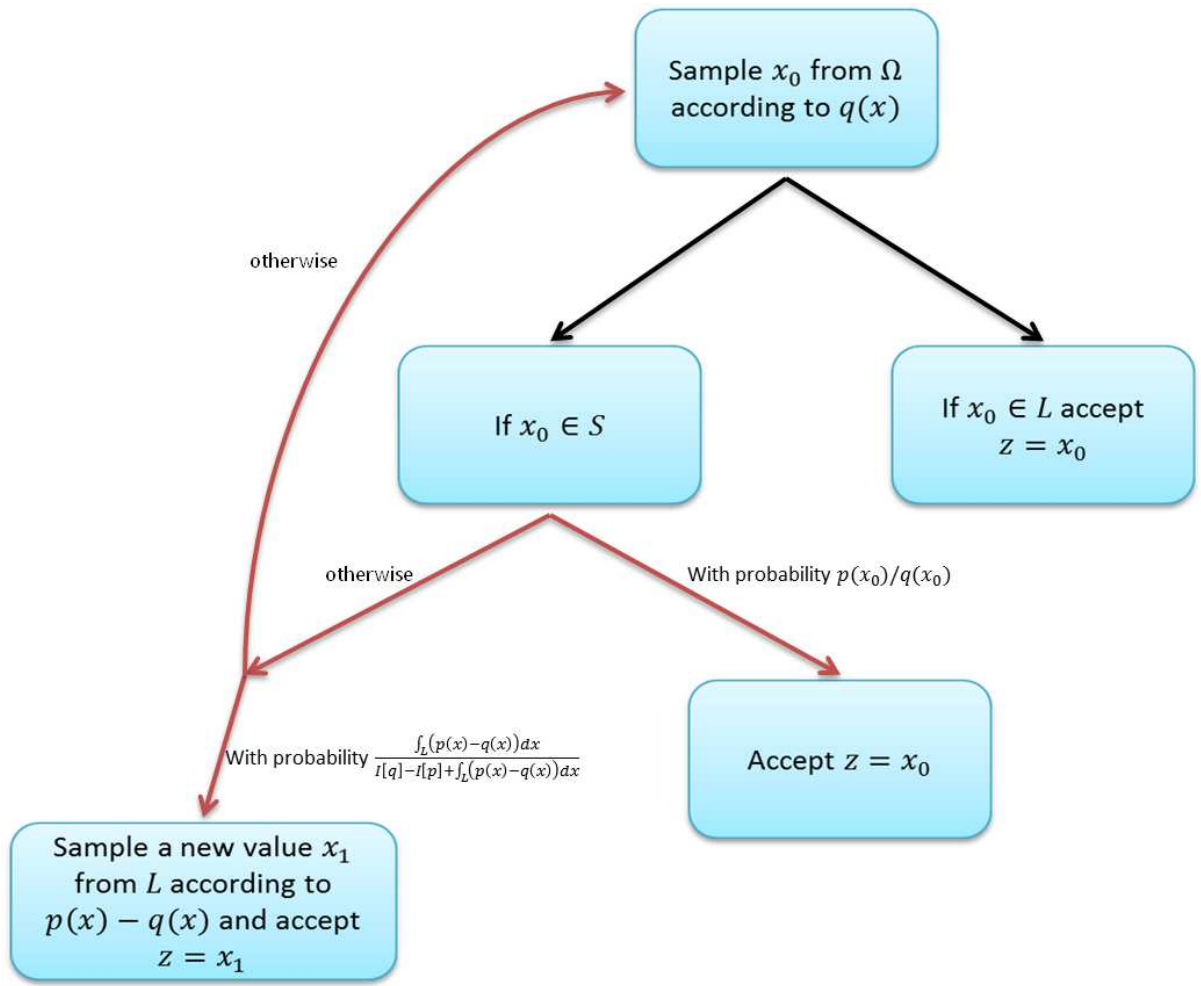} 
    \caption{Flow chart of Algorithm II of the Reduced Rejection sampling method attributing to the case $I[p] < I[q]$.}
    \label{figure:AlgorithmII}
\end{figure}

\pagebreak


\begin{thebibliography}{130}
\bibliographystyle{plain}

\bibitem{BIB:Bird1994}
G. A. Bird, {\em Molecular Gas Dynamics and the Direct Simulation of Gas Flows}, Claredon, Oxford (1994)

\bibitem{BIB:Bortz1975}
A.B. \ Bortz, M.H. \ Kalos, and J.L. \ Lebowitz. ``A  new algorithm for {Monte Carlo } simulation of {Ising} spin systems'', \textit{J. Comp. Phys.} 17 (1975), 10-18.


\bibitem{BIB:Caflisch} R.E.\ Caflisch. ``Monte Carlo and Quasi-Monte Carlo Methods'', \textit{Acta Numerica} (1998), 1--49.

\bibitem{BIB:Cao} Y.\ Cao, H.\ Li, L.R.\ Petzold, ``Efficient formulation of the stochastic simulation algorithm for chemically reacting systems", \textit{J. Chem. Phys.}, \textbf{121}, (2004), 4059--4067.

\bibitem{BIB:Deak}  I. Deak. ``An Economical Method for Random Number Generation and a Normal Generator'', \textit{Computing} 27 (1981), 113--121.

\bibitem{BIB:Gibson} M.A.\ Gibson, J.\ Bruck, `` Exact stochastic simulation of chemical systems with many species and many channels", \textit{J. Phys. Chem.}, \textbf{105}, (2000), 1876--1889.

\bibitem{BIB:GilksBestTan} W. R.\ Gilks, N. G.\ Best, and K. K. C.\ Tan, ``Adaptive rejection Metropolis sampling'', \textit{Applied Statistics}, \textbf{44}, (1995), 455--472.

\bibitem{BIB:GilksWild} W. R.\ Gilks and P.\ Wild, ``Adaptive rejection sampling for Gibbs sampling'', \textit{Applied Statistics}, \textbf{41}, (1992), 337--348.

\bibitem{BIB:Gillespie07} D.T.\ Gillespie, ``Stochastic Simulation of Chemical Kinetics'', \textit{Annu. Rev. Phys. Chem.}, \textbf{58}, (2007), 35--55.

\bibitem{BIB:Gillespie76} D.T.\ Gillespie, ``A general method for numerically simulating the stochastic time evolution of coupled chemical reactions", \textit{J. Comput. Phys.}, \textbf{22}, (1976), 403--434.

\bibitem{BIB:TableMethod} G.\ Marsaglia, W. W.\ Tsang, and J. Wang, ``Fast Generation of Discrete Random Variables'', \textit{Journal of Statistical Software}, \textbf{11}, 3 (2004).

\bibitem{BIB:Ziggurat} G.\ Marsaglia and W. W.\ Tsang, ``A fast, easily implemented method for sampling from decreasing or symmetric unimodal density functions'', \textit{SIAM Journ. Scient. and Statis. Computing}, \textbf{5}, 2 (1984), 349--359.  

\bibitem{BIB:Xorshift} G.\ Marsaglia, ``Xorshift RNGs'', \textit{Journal of Statistical Software}, \textbf{8}, 14 (2003).

\bibitem{BIB:McCollum} J.M.\ McCollum, G.D.\ Peterson, C.D.\ Cox, M.L.\ Simpson, N.F.\ Samatova, ``The sorting direct method for stochastic simulation of biochemical systems with varying reaction execution behavior", \textit{Comput. Bio. Chem.}, \textbf{30}, (2006), 39--49.

\bibitem{BIB:Oxenius} J.T. \ Oxenius, \textit{Kinetic Theory of Particles and Photons} (1986) Springer-Verlag, Berlin.

%\bibitem{BIB:UniformRatio} 	E.\ Stadlober, ``Ratio of uniforms as a convenient method for sampling from classical discrete distributions'', \textit{Proceedings of the 21st ACM conference on Winter simulation}, (1989), 484--489. 	

\bibitem{BIB:Ramaswamy09} R.\ Ramaswamy, N.\ Gonzalez-Segredo, I. F.\ Sbalzarini, ``A new class of highly efficient exact stochastic simulation algorithms for chemical reaction networks'', \textit{J. Chem. Phys.}, \textbf{130}, 244104, (2009).

\bibitem{BIB:Ramaswamy10}  R.\ Ramaswamy, I. F.\ Sbalzarini, ``A partial-propensity variant of the composition-rejection stochastic simulation algorithm for chemical reaction networks", \textit{J. Chem. Phys.}, \textbf{132}, 044102, (2010).

\bibitem{BIB:Ramaswamy11} R.\ Ramaswamy, I. F.\ Sbalzarini, ``A partial-propensity formulation of the stochastic simulation algorithm for chemical reaction networks with delays", \textit{J. Chem. Phys.}, \textbf{134}, 014106, (2011).

\bibitem{BIB:AliasMethod} M. D.\ Vose, ``A Linear Algorithm For Generating Random Numbers With a Given Distribution'', \textit{IEEE Transaction and Software Engineering}, \textbf{17}, 9 (1991), 972--975.

\bibitem{BIB:Walker} A. J. Walker. ``An efficient method for generating discrete random variables with general distributions'', \textit{ACM TOMS} 3  (1977), 253--256.  

\bibitem{BIB:Zeldovich} Y.B. Zeldovich and Y.P. \ Raizer, \textit{Physics of Shock Waves and High-Temperature Hydrodynamic Phenomena} (2002) Dover, Mineola, NY.








\end{thebibliography}
\end{document}